\documentclass[10pt]{article}

\font\smallit=cmti10

\usepackage{amssymb,latexsym,amsmath,epsfig,amsthm}
\usepackage{amsfonts}
\usepackage[justification=centering]{caption}
\usepackage{float}
\usepackage{afterpage}
\usepackage{enumerate}
\usepackage{hyperref}
\usepackage{xcolor}
\usepackage{comment}
\usepackage{mathtools}

\makeatletter

\renewcommand\section{\@startsection {section}{1}{\z@}
{-30pt \@plus -1ex \@minus -.2ex}
{2.3ex \@plus.2ex}
{\normalfont\normalsize\bfseries\boldmath}}

\renewcommand\subsection{\@startsection{subsection}{2}{\z@}
{-3.25ex\@plus -1ex \@minus -.2ex}
{1.5ex \@plus .2ex}
{\normalfont\normalsize\bfseries\boldmath}}

\renewcommand{\@seccntformat}[1]{\csname the#1\endcsname. }

\makeatother

\newtheorem{theorem}{Theorem}

\newtheorem{corollary}{Corollary}
\newtheorem{definition}{Definition}

\allowdisplaybreaks
\begin{document}

\begin{center}
{\Large \bf $d$-Translated Unit Sensitive Primes}
\vskip 20pt
{\bf Thomas Luckner}\\
{\smallit Dept.~Mathematics and Technology, 
Flagler College,}\\
{\smallit Saint Augustine, FL 32080, USA}\\
{\tt tluckner@flagler.edu}\\ 
\vskip 20pt
{\bf R. James Philpott}\\
{\smallit Dept.~Mathematics and Technology, 
Flagler College,}\\
{\smallit Saint Augustine, FL 32080, USA}\\
{\tt philpottrjames@gmail.com}\\ 
\end{center}

\vskip 5pt
\centerline{\smallit Received: , Revised: , Accepted: , Published: } 

\vskip 5pt 

\vskip 15pt

\centerline{\bf Abstract}
\vskip 5pt\noindent
A nonnegative integer $n$ is $d$-translated unit sensitive when appending $d$ zeros to the right of the number and changing the unit digit to any possible unit digit (so long as the resulting number is not $n$) results in a composite number. We find an arithmetic progression of nonnegative integers that are $d$-translated unit sensitive for any nonnegative integer $d$. We construct the arithmetic progression such that there exists $k$ consecutive primes in the progression for any positive integer $k$. The first known prime that is $d$-translated unit sensitive for any nonnegative integer $d$ is found as a consequence. We also refine the arithmetic progression such that the integers in the progression are $d$-translated unit sensitive for any nonnegative integer $d$ and are Brier numbers. The refined arithmetic progression contains $k$ consecutive primes for any choice of positive integer $k$.

\vskip 10pt
\centerline{\textbf{AMS Subject Classifications:} 11B25 (11A07, 11A51, 11A63)}

\pagestyle{myheadings} 
\thispagestyle{empty} 
\baselineskip=12.875pt 
\vskip 30pt
\section{Introduction}
In 1950, Paul Erd\H{o}s introduced the technique of covering systems as a constructive method to answer questions dealing with infinitely many integers (see \cite{erdos}). Covering systems continue to be used as a modern proof technique (see \cite{jacob}, \cite{ftj}, and \cite{jj}). Simultaneously, several open problems emerged relating to covering systems (see \cite{minmod}, \cite{maria}, and \cite{hough}). This paper uses covering systems similarly as used by Erd\H{o}s in \cite{erdos}. Thus, we first define a covering system.

\begin{definition}[Covering System] A \textit{covering system (or covering)} is a finite set of congruences,

\[\lbrace x \equiv a_{1} {\hskip -4pt}\pmod{b_{1}}, \quad x \equiv a_{2} {\hskip -4pt}\pmod{b_{2}}, \quad \ldots, \quad x \equiv a_{n} {\hskip -4pt}\pmod{b_{n}}\rbrace, \]
where $n$ is a positive integer, each $a_i$ is an integer, and each $b_i$ is a positive integer such that every integer satisfies at least one congruence in the set.
\end{definition}

In \cite{poli}, A. De Polignac conjectured all odd integers greater than 3 can be represented as the sum of a power of two and a prime. N. Romonoff showed that a positive proportion of odd numbers can be represented this way (see \cite{rom}). Let $x$ satisfy the conditions of the conjecture. That is, $x=2^k+p$ for some positive integer $k$ and prime $p$. However, several counterexamples disproved the conjecture, and consequently, Erd\H{o}s (in \cite{erdos}) constructed an arithmetic progression of odd integers $x$ where $x-2^k$ is composite for any choice of positive integer $k$. His construction utilized the covering system in first column of Table \ref{ecov}.

\begin{table}[!hbt]
    \centering
    \caption{Erd\H{o}s Covering System}
    \label{ecov}
\begin{tabular}{|c|c|}
\hline Congruence Classes for $k$  & Congruence Classes for $x$ \\
\hline \hline
$0\pmod{2}$ & $1 \pmod{3}$ \\
\hline
$0\pmod{3}$ & $1 \pmod{7}$\\
\hline
$1\pmod{4}$ & $2 \pmod{5}$\\
\hline
$3\pmod{8}$ & $8 \pmod{17}$\\
\hline
$7\pmod{12}$ & $11 \pmod{13}$\\
\hline
$23\pmod{24}$ & $121 \pmod{241}$\\
\hline
\end{tabular}
\end{table}

\noindent For each congruence class in the covering, Erd\"os ensured $x-2^k$ was composite by fixing $x$ according to the second column of Table \ref{ecov}. For example, consider the second row of Table \ref{ecov} where $k\equiv 0 \pmod{3}$ and $x\equiv 1 \pmod{7}$. Then, for the nonnegative integer $n$, 
\[
x-2^k=x-2^{0+3n}=x-2^0\cdot 2^{3n}\equiv 1-1\cdot 1\equiv 0\pmod{7}.
\] 
\noindent Thus, $7$ divides $x-2^k$ when $k\equiv 0 \pmod 3$ and $x\equiv 1 \pmod 7$. Therefore, $x-2^k$ is composite for the second row of Table \ref{ecov} so long as $x-2^k$ is not $7$.  Note that the moduli in the congruence classes for $x$ in Table \ref{ecov} are unique primes. The Chinese Remainder Theorem ensures $x-2^k$ is composite for any positive integer $x\equiv  7629217 \pmod{11184810}$ such that $x-2^k$ is not any of the prime moduli. We can express $x\equiv  7629217 \pmod{11184810}$ as the arithmetic progression
\[
x=11184810 \,m+7629217
\]
where $m$ is any nonnegative integer. 

Similarly to Erd\"os, W. Sierp\'inski used covering systems to prove that there are infinitely many positive odd integers $k$ such that $k\cdot 2^n+1$ is composite for all nonnegative integers $n$. Sierp\'inski found all integers $k$ such that 
\[
k=36893488147419103230\,m+ 15511380746462593381 
\] 
for some nonnegative integer $m$ have this property (see \cite{sier}). Consequently, we refer to integers of this form as \textit{Sierp\'inski numbers}. Likewise, H. Riesel proved that there are infinitely many positive odd integers $k$ such that $k\cdot 2^n-1$ is composite for all nonnegative integers $n$. Riesel found all integers $k$ such that
\[
k=11184810\,m+509203
\] 
for some nonnegative integer $m$ have this property (see \cite{ries}). Consequently, integers of this form are referred to as \textit{Riesel numbers}. Clavier found the smallest known example of a \textit{Brier number} (a number that is both a Sierp\'inski number and a Riesel number).  By way of methods similar to Erd\"os, Sierp\'inski, and Riesel, Clavier found every integer $k$ such that
\[
k=3316923598096294713661 + 3770214739596601257962594704110\, m 
\]
for some nonnegative integer $m$ is a Brier number (see \cite{slobrier} and \cite{jacob}). In \cite{jacob}, Filaseta and Juillerat conclude that the arithmetic progression contains infinitely many prime Brier numbers. For other recent results pertaining to Sierp\'inski, Riesel, and Brier numbers see: \cite{barn}, \cite{harr1}, \cite{brucald}, \cite{fil}, \cite{ftj}, \cite{harr2}, \cite{slosier}, and \cite{slories}. 

While our main result pertains to Brier numbers, there is another class of integers fundamental to the results of the paper. First, we define the following class of integers.

\begin{definition}[Unit Sensitive]
Let any positive integer be unit sensitive if changing the unit digit to any other base 10 digit results in a composite number. 
\end{definition}
We show that the integer $97$ is unit sensitive. Changing the unit digit results in one of the following integers:
\[
90, 91, 92, 93, 94, 95, 96, 98, 99.
\]
 All of the above integers are composite. Thus, since $97$ is prime, $97$ is a unit sensitive prime. 
This is not a particularly rare occurrence. In fact, we  check if a prime is unit sensitive by using only nine primality checks. We extend the notion of unit sensitivity in the following way.
\begin{definition}[$d$-Translated Unit Sensitive]
    For any nonnegative integer $d$, a positive integer, $n$, is $d$-translated unit sensitive if appending $d$ zeros to the right of the number and changing the unit digit to any possible unit digit results in a composite number (so long as the resulting number is not $n$).
    
    Note that, in the case $d=0$, $n$ is unit sensitive.
\end{definition}
To illustrate this construction, we show $97$ is not $1$-translated unit sensitive. We need only check if the following integers are composite:
\[
970, 971, 972, 973, 974, 975, 976, 977, 978, 979.
\]
Since 971 and 977 are prime, we say $97$ is not $1$-translated unit sensitive. 
By similar argument, the first $1$-translated unit sensitive prime is $53$. By repeating this procedure for $d=2$ we can show that $89$ is the first prime that is both $1$-translated unit sensitive and $2$-translated unit sensitive. For simplicity, when an integer, $n$, is $d$-translated unit sensitive for all integers $d \in [a,b]$ for some integers $a$ and $b$, we say $n$ is $[a,b]$-translated unit sensitive. Thus, $89$ is the first $[1,2]$-translated unit sensitive prime. In fact, as seen in Table \ref{comp1} and Table \ref{comp2} respectively, none of the first $10^6$ primes are $[1,4381]$-translated unit sensitive and none of the first $10^6$ primes are $[0, 3191]$-translate unit sensitive.

\begin{table}[!hbt]
    \centering
    \caption{Largest integer $b$ such that, of the first $N$ primes, a prime is $[1,b]$-translated unit sensitive}
    \label{comp1}
    \begin{tabular}{|c|c|c|}
    \hline $N$ & $b$ & $[1,b]$-translated unit sensitive prime\\
    \hline \hline
     $10^2$ & $9$ & $443$\\
     \hline
     $10^3$ & $38$ & $7829$\\
     \hline
     $10^4$ & $269$ & $52121$\\
     \hline
     $10^5$ & $513$ & $390893$\\
     \hline
     $10^6$ & $4380$ & $2049797$\\
     \hline
    \end{tabular}
\end{table}

\begin{table}[!hbt]
    \centering
 \caption{Largest integer $b$ such that, of the first $N$ primes, a prime is $[0,b]$-translated unit sensitive}
    \label{comp2}
    \begin{tabular}{|c|c|c|}
   \hline $N$ & $b$ & $[0,b]$-translated unit sensitive prime\\
    \hline \hline
     $10^2$ & $5$ & $113$\\
     \hline
     $10^3$ & $11$ & $5693$\\
     \hline
     $10^4$ & $51$ & $37379$\\
     \hline
     $10^5$ & $163$ & $951161$\\
     \hline
     $10^6$ & $3190$ & $10721819$\\
     \hline
    \end{tabular}
\end{table}

Similar to the unit sensitive condition, the $d$-translated unit sensitive condition requires only ten primality checks for each choice of $d\ge 1$. Thus, a computer checks whether an integer is $d$-translated unit sensitive for a finite set of $d$ in a reasonable amount of time. To refine this notion further such that computer computation cannot reasonably check the condition, we define an integer $n$ to be  $[0,\infty)$-translated unit sensitive when $n$ is $d$-translated unit sensitive for any nonnegative integer $d$. Since $b$ is increasing exponentially, Table \ref{comp2} suggests there may not be any $[0,\infty)$-translated unit sensitive primes. Contrary to this observation, we prove the following results and corollaries.

\begin{theorem}\label{main1}
There are infinitely many primes that are $d$-translated unit sensitive for any nonnegative integer $d$. In fact, there are infinitely many primes of the form, $Am+B$ for some nonnegative integer $m$, that are $[0,\infty)$-translated unit sensitive for
\[
 A=124703827160493827160493827160493827160493827160493827035790 \text{ and }
\]
\[
B=41459060189171787548442999328384678040412832671445258454633.
\]
\end{theorem}
Since $B$ is prime, $B$ is the first known example of a $[0,\infty)$-translated unit sensitive prime. Theorem \ref{main2} adds the property of Brier numbers.

\begin{theorem}\label{main2}
There are infinitely many primes that are $d$-translated unit sensitive for any nonnegative integer $d$ (or $[0,\infty)$-translated unit sensitive) and are also Brier numbers. In fact, there are infinitely many primes of the form, $Am+B$ for some nonnegative integer $m$, that are $[0,\infty)$-translated unit sensitive and Brier numbers for

\[
A=105882605432768363742137390170084901366252883723703759364986698321995
\]
\[
2824339882555233842344362162520787653763315998320477054896263434901090 \text{ and }
\]
\[
B=292134647723778159914962540207823140335126344033210264225285566059959
\]
\[
956372273970360496308095920463669048323126517085888353814417281244783.
\]
\end{theorem}
The first known prime that is both $[0,\infty)$-translated unit sensitive and a Brier number is $A (11)+B$, that is, 
\[
119392212453282981715500754589171622906229435536406237943738223814794410
\]
\[
24110982077932762096079708192333239719602498611135957673315065156773.
\]
The following corollaries are consequences of the construction necessary for every integer in the arithmetic progression in Theorem \ref{main2} to be a Brier numbers.

\begin{corollary}\label{cor1}
The arithmetic progression, $Am+B$ for any nonnegative integer $m$, contains infinitely many prime Sierp\'inski numbers where
\[
36893488147419103230\,m + 10691053702625738573.
\]
\end{corollary}
\begin{corollary}\label{cor2}
The arithmetic progression, $Am+B$ for any nonnegative integer $m$, contains infinitely many prime Riesel numbers where
\[
A=3189762509404991434026267111935081578564848099075598418737381122 \text{ and }
\]
\[
B=1464481329052452484012410951111069258728979724777952402192519099.
\]
\end{corollary}
\begin{corollary}\label{cor3}
The arithmetic progression, $Am+B$ for any nonnegative integer $m$, contains infinitely many prime Brier numbers where
\[
A=19613577555635811210751525911492599150619172808998516420233576073634
\]
\[
048313995204010 \text{ and }
\]
\[
B=76099738321991043329005201186973924664733439211950059278085420869117
\]
\[
02037440072483.
\]
\end{corollary}
Corollaries \ref{cor1}, \ref{cor2}, and \ref{cor3} provide arithmetic progressions of Sierp\'inski, Riesel and Brier numbers respectively different from those mentioned in \cite{jacob}, \cite{ries}, and \cite{sier}. All of the above results extend using a result of D. Shiu (see \cite{shiu}) or a stronger result of J. Maynard (see \cite{maynard}) just as was done in \cite{ftj}. We demonstrate these extensions using Theorem \ref{main2}.

\begin{corollary}[Shiu Corollary]\label{shiucor}
For every positive integer $k$, there exist $k$ consecutive primes each of which is a $[0, \infty)$-translated unit sensitive prime and a Brier number.
\end{corollary}
\begin{corollary}[Maynard Corollary]\label{maycor}
For every positive integer $k$, there exist $k$ consecutive primes $p_{\ell}$, $p_{\ell+1}$,$\ldots$, $p_{\ell+k-1}$, each of which is both a $[0, \infty)$-translated unit sensitive prime and a Brier number. 
\end{corollary}

\section{Construction}
To prove Theorems \ref{main1} and \ref{main2}, we construct an arithmetic progression $Am+B$ that satisfies two main properties: $Am+B$ is $[0,\infty)$- translated unit sensitive and $Am+B$ is a Brier number. 

To begin the construction, fix $A$ and $B$ as positive integers. In order for the arithmetic progression to contain infinitely many primes, fix $A$ and $B$ to be relatively prime (Dirichlet). In addition to these conditions, $A$ and $B$ will be constructed such that the following criteria are met.\\

\noindent \textbf{Criterion (1):} For all integers $n\in[0,9]$, each integer in the set
    \[
    \mathcal{D}_n=\mathcal{D}_n(A, B)=\lbrace (Am+B)\cdot 10^d +n \text{ : } m\in \mathbb{Z}^+\cup \lbrace 0\rbrace \text{ and } d\in\mathbb{Z}^+  \rbrace
    \]
    is divisible by at least one prime, $p$, such that $p$ divides $A$.\\

\noindent \textbf{Criterion (2):} For all integers $x\in\{-3, -2, -1, 1, 2, 3, 4, 5, 6\}$, each integer in the set
    \[
    \mathcal{D}_x=\mathcal{D}_x(A, B)=\lbrace (Am+B)+x \text{ : } m\in \mathbb{Z}^+\cup \lbrace 0\rbrace  \rbrace
    \]
    is divisible by at least one prime, $p$, such that $p$ divides $A$.\\
    
\noindent \textbf{Criterion (3):} Each integer in the set
    \[
    \mathcal{S}=\mathcal{S}(A, B)=\lbrace (Am+B)\cdot 2^k +1 \text{ : } m,k\in \mathbb{Z}^+\cup \lbrace 0\rbrace   \rbrace
    \]
    is divisible by at least one prime, $p$, such that $p$ divides $A$.\\
    
\noindent \textbf{Criterion (4):} Each integer in the set
    \[
    \mathcal{R}=\mathcal{R}(A, B)=\lbrace (Am+B)\cdot 2^k -1 \text{ : } m,k\in \mathbb{Z}^+\cup \lbrace 0\rbrace   \rbrace
    \]
    is divisible by at least one prime, $p$, such that $p$ divides $A$.\\

First, note that Criterion (1) ensures $Am+B$ is $[1,\infty]$-translated unit sensitive. Similarly Criterion (3) and (4) ensure $Am+B$ is a Sierp\'inski number and Riesel number respectively when $Am+B$ is odd. Thus, the combination of Criteria (3) and (4) ensures $Am+B$ is a Brier number when $Am+B$ is odd. 

Lastly, note that Criterion (2) ensures that $Am+B$ is unit sensitive if $Am+B\equiv 3 \pmod {10}$. In our construction we ensure $Am+B\equiv 1 \pmod 2$ and $Am+B \equiv 3 \pmod 5$, or equivalently $Am+B\equiv 3 \pmod {10}$ (see Appendix Table \ref{d=0}). Assuming $Am+B\equiv 3 \pmod {10}$, we know $Am+B+x$ represents all the possible unit digits changes if $x$ is any nonzero integer in $[-3,6]$ as Table \ref{crit2} indicates.

\begin{table}[!hbt]
    \centering
    \caption{$\left(Am+B\right)+x$ under the conditions in Criterion (2)}
    \label{crit2}
    \begin{tabular}{|c|c|}
    \hline $x$ & $\left(Am+B\right)+x \pmod {10}$\\
    \hline \hline
     $-3$ &  $0 \pmod {10}$\\
     \hline
      $-2$ &  $1 \pmod {10}$\\
     \hline
      $-1$ &  $2 \pmod {10}$\\
     \hline
      $1$ &  $4 \pmod {10}$\\
     \hline
      $2$ &  $5 \pmod {10}$\\
     \hline
      $3$ &  $6 \pmod {10}$\\
     \hline
      $4$ &  $7 \pmod {10}$\\
     \hline
      $5$ &  $8 \pmod {10}$\\
     \hline
      $6$ &  $9 \pmod {10}$\\
     \hline
   
    \end{tabular}
\end{table}
 Thus, Criterion (2) corresponds to $Am+B$ being unit sensitive. Therefore, Criteria (1) and (2) ensures that $Am+B$ is $[0,\infty)$-translated unit sensitive. 

It is tempting to combine Criteria (1) and (2) into a single criterion by allowing $d=0$ for Criterion (1). However, this implies $Am+B\equiv 0 \pmod{10}$ and $Am+B$ is composite for any choice of nonnegative $m$. For $Am+B$ to be prime for infinitely many nonnegative $m$ (as Theorem \ref{main1} and \ref{main2} as well as Corollary \ref{cor1} - \ref{maycor} require), we need $Am+B\equiv 1, 3, 7, \text{ or } 9 \pmod {10}$. To demonstrate this restriction, consider the $[0,1]$-translated unit sensitive prime $113$.  We know $113$ is $1$-translated unit sensitive since 
\[
1130, 1131, 1132, 1133, 1134, 1135, 1136, 1137, 1138, \text{ and } 1139
\]
are all composite. In other words, since $113\cdot 10^1+n=1130+n$ is composite for all integers $n\in [0,9]$, $113$ is $1$-translated unit sensitive. The same method applies for any positive integer $d$. When $d=0$, the construction of Criterion (1) requires
\[
113, 114, 115, 116, 117, 118, 119, 120, 121, \text{ and } 122
\]
to be composite. However, we know $113$ is unit sensitive since  
\[
110, 111, 112, 114, 115, 116, 117, 118, \text{ and } 119
\]
are all composite. In other words, since $113+x$ is composite for all nonzero integers $x\in [-3,6]$, $113$ is unit sensitive. Therefore, the case when $d=0$ must be made separate from Criterion (1) to handle the case of unit sensitivity. As mentioned previously, Criterion (2) ensures $Am+B$ is unit sensitive if $Am+B\equiv 3 \pmod {10}$ which is true given our construction of coverings. 

\section{Coverings}\label{s3}
This section addresses the construction of the coverings. The coverings are constructed such that Criteria (1) - (4) is satisfied, $A$ and $B$ positive, and $A$ and $B$ are relatively prime. 

We start with Criterion (1). Notice that $(Am+B)\cdot 10^d+n$ is even for all positive $d$ when $n$ is even. Similarly, when $n=5$, $(Am+B)\cdot 10^d+n$ is divisible by $5$. Thus, coverings are only necessary for $n\in\{1, 3, 7, 9\}$ to satisfy Criterion (1). For $n=1$ and $n=7$, Criterion (1) is satisfied when 
\[
A\equiv 0 \pmod 3 \text{ and } B\equiv 2 \pmod {3}.
\]
Since $(Am+B)\cdot 10^d+n \equiv (Am+B)+n \pmod {3}$ for any positive integer $d$, these conditions modulo 3 give us
\[
(Am+B)+n\equiv 2+n \pmod {3}.
\]
Therefore, when $n\equiv 1 \pmod 3$, we have $(Am+B)\cdot10^d+n$ is divisible by 3 and composite as Appendix Table \ref{n1n7} indicates. The remaining two cases of Criterion (1), $n=3$ and $n=9$, require more complicated coverings we discuss later in this section.

Next we address Criterion (2). Notice that when $Am+B\equiv 3 \pmod  {10}$, if $x$ is odd we have $(Am+B)+x$ is even. Also, if $x=2$, $(Am+B)+x$ is divisible by 5. Thus, coverings are only necessary for $x\in\{-2, 4, 6\}$ for Criterion (2) to be satisfied. Similarly, if $x=-2$ or $x=4$ we have
\[
(Am+B)+x\equiv 2+x\equiv 0\pmod 3
\]
when $A\equiv 0 \pmod 3$ and $B\equiv 2 \pmod 3$.
We already imposed this condition on $A$ and $B$ for Criterion (1) and we need both Criteria (1) and (2) for $Am+B$ to be $[0,\infty)$-translated unit sensitive. Thus, we know the cases of $x=-2$ and $x=4$ are satisfied for Criterion (2). Therefore, the only case remaining for Criterion (2) is $x=6$. For this remaining case, we take 
\[
A\equiv 0 \pmod 7 \text{ and } B\equiv 1 \pmod 7. 
\]
Under these conditions, we have 7 divides $(Am+B)+6$. This completes Criterion (2) as Appendix Table \ref{d=0} indicates. 

The remaining choice of $n$ to satisfy Criterion (1) completely are $n=3$ and $n=9$. We also must satisfy Criteria (3) and (4). We create a covering for each of the two remaining cases of Criterion (1), Criterion (3), and Criterion (4) as seen in Tables \ref{n3} - \ref{Scov}. It is important to note that there are four moduli that are used in multiple ways. For example, we noted that when 
\[
A\equiv 0 \pmod 3 \text{ and }B\equiv 2 \pmod 3
\]
Criterion (1) is satisfied for $n=1$ and $n=7$. Similarly, Criterion (2) was satisfied for $x=-2$ and $x=4$ under these conditions. Thus, the modulus 3 was used in four different ways.  

Assuming $A\equiv 0 \pmod M$ for each modulus $M$, Table \ref{reuse} provides the congruence class $B\equiv \beta \pmod M$ used in more than one way and lists the specific criteria where it is used. In the case of Criteria (1) and (2), the specific $n$ and $x$ are also provided. In the case that the congruence class $B\equiv \beta \pmod M$ satisfies only some choices of $d$ (as in Criterion (1)) or $k$ (as in Criteria (3) and (4)), a colon followed by the specific congruence class it satisfies. For example, in Table \ref{reuse} the column for $B\equiv 1 \pmod 7$ has Criteria (1) $n=9$: $d\equiv 5 \pmod 6$ in the first row. This means that Criteria (1) is satisfied for $n=9$ when $d\equiv 5 \pmod 6$. 

\begin{table}[!hbt]
    \centering
    \caption{All Reused Congruence Classes on $Am+B$ Part 1}
    \label{reuse}
    \begin{tabular}{|c|c|}
    \hline $B\equiv 2\pmod 3$ & $B\equiv 1\pmod 7$  \\
    \hline \hline
    Criterion (1) $n=1$ and $n=7$ & Criterion (1) $n=9$: $d\equiv 5 \pmod 6$    \\
     \hline
     Criterion (2) $x=-2$ and $x=4$ & Criterion (2) $x=6$  \\
    \hline 
  Criterion (3): $k\equiv 0 \pmod 2$ & Criterion (4): $k\equiv 0 \pmod 3$  \\
   \hline 
   Criterion (4): $k\equiv 1 \pmod 2$ & -  \\
   \hline
    \end{tabular}
\end{table}
\addtocounter{table}{-1}
\begin{table}[!hbt]
    \centering
    \caption{All Reused Congruence Classes on $Am+B$ Part 2}
    \begin{tabular}{|c|c|}
    \hline $B\equiv 3\pmod {11}$ & $B\equiv 3\pmod {10}$\\
    \hline \hline
     Criterion (1) $n=3$: $d\equiv 1 \pmod 2$  & Criterion (2)   \\
     \hline
    Criterion (4): $k\equiv 2 \pmod {10}$ & Criterion (3): $k\equiv 3 \pmod 4$\\
    \hline 
    
    \end{tabular}
\end{table}

Fix $n$ such that $n=3$ or $n=9$ to finish the construction of $A$ and $B$ to satisfy Criterion (1). Let $d\equiv a \pmod b$ be a congruence class in the covering with unique (with the exception of any primes from Table \ref{reuse}) associated prime $p$. Then $d=a+b\ell$ for some integer $\ell$. By selecting $b$ such that the order of 10 modulo $p$ divides $b$, we have  
\[
(Am+B)\cdot 10^d+n\equiv (Am+B)\cdot 10^{a+b\ell}+n \equiv (Am+B)\cdot 10^a +n \pmod {p}.
\]
Assuming $p$ divide $A$, we have $p$ divides $(Am+B)\cdot 10^d+n$ when 
\[
(Am+B)\cdot 10^a+n\equiv 0 \pmod p \Rightarrow B\equiv \dfrac{-n}{10^a} \pmod p.
\]
We have $B$ exists and is relatively prime to $A$ when $p\neq 2$ and $p\neq 5$. Formally, for each congruence class in the covering system $d\equiv a \pmod b$ and associated prime $p$ we let 
\[
A\equiv 0\pmod p \text{ and } B\equiv \dfrac{-n}{10^a} \pmod p
\]
such that $p$ divides $(Am+B)\cdot 10^d+n$. For both $n=3$ and $n=9$, this construction remains the same and completes Criterion (1). The construction of the coverings for Criterion (3) and Criterion (4) follows similarly assuming $p\neq 2$. Assuming $A\equiv 0 \pmod p$ for each prime $p$, Table \ref{constr} provides the construction of congruence classes $B$ must satisfy for the remaining four cases: Criterion (1) when $n=3$, Criterion (1) when $n=9$, Criterion (3), and Criterion (4). 

\begin{table}[!hbt]
    \centering
    \caption{Construction of Congruence Classes $B$ Must Satisfy}
    \label{constr}
    \begin{tabular}{|c|c|c|}
    \hline Remaining Cases & Congruence Class in Covering & $B\equiv \beta \pmod p$ \\
    \hline \hline
    Criterion (1) $n=3$ & $d\equiv a \pmod b$ & $B\equiv \left(-3/10^a\right)\pmod p$ \\
     \hline
     Criterion (1) $n=9$ & $d\equiv a \pmod b$ & $B\equiv \left(-9/10^a\right)\pmod p$  \\
     \hline
     Criterion (3) & $k\equiv a \pmod b$ & $B\equiv \left(-1/2^a\right)\pmod p$   \\
     \hline
     Criterion (4) & $k\equiv a \pmod b$& $B\equiv \left(1/2^a\right)\pmod p$  \\
     \hline
    \end{tabular}
\end{table}

Since all primes used are unique except for those in Table \ref{reuse}, the construction of the four remaining cases is complete.  Please see Appendix Tables \ref{d=0}--\ref{Rcov} for a list of the coverings that complete Criteria (1) - (4). 

Now we explicitly construct $A$ and $B$ according to Appendix Tables \ref{d=0}--\ref{Rcov}.  We begin with Criterion (1). Let
\[
\mathcal{P}_{\mathcal{D}_1}=\mathcal{P}_{\mathcal{D}_7}=\{3 \} \text{ and } \mathcal{B}_{\mathcal{D}_1}=\mathcal{B}_{\mathcal{D}_7}=\{2\pmod 3\}.
\] 
For $n=3$ and $n=9$ in Criterion (1), using the construction defined in Table \ref{constr}, choose
\[
\mathcal{P}_{\mathcal{D}_3}=\{p \text{ : } p \text{ is a prime in column two of Appendix Table \ref{n3}} \},
\] 
\[
\mathcal{P}_{\mathcal{D}_9}=\{p \text{ : } p \text{ is a prime in column two of Appendix Table \ref{n9}} \},
\] 
\[
\mathcal{B}_{\mathcal{D}_3}=\{ \beta \pmod p \text{ : } \beta \pmod p \text{ is in column three of Appendix Table \ref{n3}}\}, \text{ and}
\]
\[
\mathcal{B}_{\mathcal{D}_9}=\{ \beta \pmod p \text{ : } \beta \pmod p \text{ is in column three of Appendix Table \ref{n9}}\}, \text{ and}
\]
To satisfy Criterion (2), choose  
\[
\mathcal{P}_{\mathcal{D}_{-2}}=\{2, \text{ } 3, \text{ } 5\}, \text{ }\mathcal{P}_{\mathcal{D}_{4}}=\{2, \text{ } 3, \text{ } 5\},  \text{ and }
\mathcal{P}_{\mathcal{D}_{6}}=\{2, \text{ } 5, \text{ } 7\}.
\]
Similarly, choose
\[
\mathcal{B}_{\mathcal{D}_{-2}}=\{1 \pmod 2,\text{ } 2\pmod 3, \text{ } 3\pmod 5\},
\]
\[
\mathcal{B}_{\mathcal{D}_{4}}=\{1 \pmod 2, \text{ } 2\pmod 3, \text{ } 3\pmod 5\}, \text{ and } 
\]
\[
\mathcal{B}_{\mathcal{D}_{-2}}=\{1 \pmod 2, \text{ } 3\pmod 5, \text{ } 1\pmod 7\}.
\]
To satisfy  Criteria (3) and (4), we use the construction defined in Table \ref{constr} and the coverings in Appendix Tables \ref{Scov} and \ref{Rcov}. Choose
\[
\mathcal{P}_{\mathcal{S}}=\{p \text{ : } p \text{ is a prime in column two of Appendix Table \ref{Scov}} \},
\] 
\[
\mathcal{P}_{\mathcal{R}}=\{p \text{ : } p \text{ is a prime in column two of Appendix Table \ref{Rcov}} \}, 
\]
\[
\mathcal{B}_{\mathcal{S}}=\{ \beta \pmod p \text{ : } \beta \pmod p \text{ is in column three of Appendix Table \ref{Scov}}\}, \text{ and }
\] 
\[
\mathcal{B}_{\mathcal{R}}=\{ \beta \pmod p \text{ : } \beta \pmod p \text{ is in column three of Appendix Table \ref{Rcov}}\}.
\]
We find $A$ and $B$ such that Criteria (1) - (4) are satisfied by defining 
\[
\mathcal{P}=\left(\displaystyle\bigcup_{n\in\{1, 3, 7, 9\}} \mathcal{P}_{\mathcal{D}_n}\right) \cup \left( \displaystyle\bigcup_{x\in\{-2, 4, 6\}} \mathcal{P}_{\mathcal{D}_x}\right)\cup \mathcal{P}_{\mathcal{S}} \cup \mathcal{P}_{\mathcal{R}}
\]
and
\[
\mathcal{B}=\left(\displaystyle\bigcup_{\text{all }n \text{ in  }\{1, 3, 7, 9\}} \mathcal{B}_{\mathcal{D}_n}\right) \cup \left( \displaystyle\bigcup_{\text{all }x \text{ in  }\{-2, 4, 6\}} \mathcal{B}_{\mathcal{D}_x}\right)\cup \mathcal{B}_{\mathcal{S}} \cup \mathcal{B}_{\mathcal{R}}.
\]
As mentioned previously, we choose $A\equiv 0 \pmod p$ for all $p\in \mathcal{P}$. Since each prime is unique, we choose
\[
A=\prod_{p\in\mathcal{P}}p.
\]
We choose $B$ such that $B\equiv \beta \pmod p$ for all $\beta\pmod p \in \mathcal{B}$. By the Chinese Remainder Theorem, given that each prime $p$ is unique, there exists some positive integer $\beta$ such that 
\[
B\equiv  \beta_{T2} \pmod A.
\]
We choose $B=\beta_{T2}$. By construction $A$ and $B$ are relatively prime. Thus,  Dirichlet's Theorem proves Theorem \ref{main2}.  To prove Theorem \ref{main1}, define 

\[
\mathcal{P}=\left(\displaystyle\bigcup_{n\in\{1, 3, 7, 9\}} \mathcal{P}_{\mathcal{D}_n}\right) \cup \left( \displaystyle\bigcup_{x\in\{-2, 4, 6\}} \mathcal{P}_{\mathcal{D}_x}\right)
\]
and
\[
\mathcal{B}=\left(\displaystyle\bigcup_{\text{all }n \text{ in  }\{1, 3, 7, 9\}} \mathcal{B}_{\mathcal{D}_n}\right) \cup \left( \displaystyle\bigcup_{\text{all }x \text{ in  }\{-2, 4, 6\}} \mathcal{B}_{\mathcal{D}_x}\right).
\]
Once again, we choose 
\[
A=\prod_{p\in\mathcal{P}}p
\]
and $B$ such that $B\equiv \beta \pmod p$ for all $\beta \pmod p \in \mathcal{B}$. The Chinese Remainder Theorem ensures there exists an integer $\beta_{T1}$ such that 

\[
B\equiv \beta_{T1} \pmod A.
\]
We choose $B=\beta_{T1}$. Thus, Dirichlet's Theorem proves Theorem \ref{main1}.  Corollary \ref{cor1} is a consequence of using 
\[
\mathcal{P}=\mathcal{P}_\mathcal{S} \text{ and } \mathcal{B}=\mathcal{B}_\mathcal{S}
\]
to construct $A$ and $B$ in the same way. Corollary \ref{cor2} is a consequence of using 
\[
\mathcal{P}=\mathcal{P}_\mathcal{R} \text{ and } \mathcal{B}=\mathcal{B}_\mathcal{R}
\]
to construct $A$ and $B$ in the same way.  Lastly, Corollary \ref{cor3} is a consequence of using 
\[
\mathcal{P}=\mathcal{P}_\mathcal{S}\cup\mathcal{P}_\mathcal{R} \text{ and } \mathcal{B}=\mathcal{B}_\mathcal{S}\cup\mathcal{B}_\mathcal{R}
\]
to construct $A$ and $B$ in the same way.  

\section{Note on Computational Methods}
To find the primes $p$ used as moduli in $B\equiv \beta \pmod p$ for Criterion (1) when $n=3$ and Criterion (1) when $n=9$, we used the computer program Magma V2.25-7 to factor $\Phi_{b}(10)$ into prime factors. Given the congruence class in the covering, say $B\equiv \beta \pmod p$ relating to the congruence class $d\equiv a \pmod b$, the prime $p$ is a factor of $\Phi_{b}(10)$ . For Criteria (3) and (4) we used the same method to factor $\Phi_{b}(2)$ into prime factors for $k\equiv a \pmod b$.

We confirm a covering by ensuring every integer in $[0,\ell-1]$ belongs to at least one of the congruence classes in the covering when $\ell$ is the least common multiple of all the moduli $b$. For example, in column one of Appendix Table \ref{n3}, the congruence classes $d\equiv 1 \pmod 2$, $d\equiv 2\pmod 4$, $d\equiv 4\pmod 8$, and $d\equiv 0 \pmod 8$ form a covering. We confirm this by checking all integers in $[0,7]$ satisfy at least one of these congruence classes. For this paper, utilizing the computer program Maple 2019, this method is sufficient.  For larger coverings, this method can become an issue. For an alternate method of covering verification, see \cite{ftj}.

\newpage
\vskip 10pt
\centerline{\textbf{\large Appendix}}
\vskip 7pt \noindent

To aid in verifying the computations in this paper, all the data in this appendix is found in \cite{data}  suitable for computations. 
This appendix begins with Table~\ref{d=0} which illustrates which primes are used in satisfying Criteria (2). Table \ref{d=0} also includes $\beta \pmod p$ such that we choose $B\equiv \beta \pmod p$ to satisfy Criterion (2). 

Tables~\ref{n1n7}--\ref{n9} contain similar information such that Criterion (1) is satisfied when $n=1$ and $n=7$, $n=3$, and $n=9$ respectively. The combination of Tables~\ref{n1n7}--\ref{n9} completes Criteria (1). The combination of Tables~\ref{d=0}--\ref{n9} is used to find $A$ and $B$ such that $Am+B$ is $[0,\infty)$-translated unit sensitive for any nonnegative integer $m$.

Tables~\ref{Scov} and \ref{Rcov} contain similar information such that Criteria (3) and (4) are satisfied respectively. The combination of Tables~\ref{Scov} and \ref{Rcov} completes Criteria (3) and (4) and is used to find $A$ and $B$ such that $Am+B$ is a Brier number for any nonnegative integer $m$.

There are several entries in the tables of this appendix that are notated with an asterisk (*). This is used to note a prime $p$ and congruence class $\beta \pmod p$ such that $B\beta \pmod p$ that was used in multiple coverings. For a list of all primes and congruence classes used in multiple coverings see Table \ref{reuse} in Section \ref{s3}. A prime was reusable in another covering so long as the congruence class $\beta \pmod p$ remains the same. For example in Table~\ref{n1n7} the prime 3 is used so that $B\equiv 2 \pmod{3}$. This prime and congruence class are used again in both Table~\ref{Scov} and Table~\ref{Rcov} with $B\equiv 2\pmod {3}$.

\newpage

\begin{center}
\begin{table}[!hbt]
    \centering
    \caption{Covering Systems for Criteria (2)}
    \label{d=0}
    \begin{tabular}{|c|c|c|}
    \hline $x$ & $p$ & $B \equiv \beta \pmod p$\\
    \hline \hline
     - & $2^*$ & $1\pmod {2}$\\
     \hline
     - & $5^*$ & $3\pmod {5}$\\
     \hline
      $x=-2$ and $x=4$ & $3^*$ & $2\pmod {3}$\\
     \hline
      $x=6$ & $7^*$ & $1\pmod {7}$\\
     \hline

    \end{tabular}
\end{table}
\end{center}

\begin{center}
\begin{table}[!hbt]
    \centering
    \caption{Covering System for $n=1$ and $n=7$}
    \label{n1n7}
    \begin{tabular}{|c|c|c|}
    \hline Covering & $p$ & $B\equiv \beta  \pmod p$\\
    \hline \hline
      $d\equiv 0\pmod 1$ & $3^*$ & $2\pmod {3}$\\
     \hline

    \end{tabular}
\end{table}
\end{center}

\begin{center}
\begin{table}[!hbt]
    \centering
    \caption{Covering System for $n=3$}
    \label{n3}
    \begin{tabular}{|c|c|c|}
    \hline Covering & $p$ & $B\equiv \beta \pmod p$\\
    \hline \hline
     $d\equiv 1 \pmod {2}$ & $11^*$ & $3\pmod {11}$\\
     \hline
     $d\equiv 2 \pmod {4}$ & 101 & $3\pmod {101}$\\
     \hline
     $d\equiv 4 \pmod {8}$ & 73 & $3\pmod {73}$\\
     \hline
     $d\equiv 0 \pmod {8}$ & 137 & $134\pmod {137}$\\
     \hline
    \end{tabular}
\end{table}
\end{center}

\begin{center}
\begin{table}[!hbt]
    \centering
    \caption{Covering System for $n=9$}
    \label{n9}
    \begin{tabular}{|c|c|c|}
    \hline Covering & $p$ & $B\equiv \beta  \pmod p$\\
    \hline \hline
     $d\equiv 1 \pmod {3}$ & 37 & $25\pmod {37}$\\
     \hline
     $d\equiv 2 \pmod {6}$ & 13 & $12\pmod {13}$\\
     \hline
     $d\equiv 5 \pmod {6}$ & $7^*$ & $1\pmod {7}$\\
     \hline
     $d\equiv 6 \pmod {9}$ & 333667 & $324667\pmod {333667}$\\
     \hline
      $d\equiv 3 \pmod {18}$ & 19 & $4\pmod {19}$\\
     \hline
      $d\equiv 12 \pmod {18}$ & 52579 & $43588\pmod {52579}$\\
     \hline
      $d\equiv 9 \pmod {27}$ & 757 & $514\pmod {757}$\\
     \hline
          $d\equiv 18 \pmod {27}$ & 440334654777631 & $440325654777631\pmod {440334654777631}$\\
     \hline
     $d\equiv 0 \pmod {54}$ & 70541929 & $70541920\pmod {70541929}$\\
     \hline
     $d\equiv 27 \pmod {54}$ & 14175966169 & $9\pmod {14175966169}$\\
     \hline
    \end{tabular}
\end{table}
\end{center}

\newpage

\begin{center}
\begin{table}[!hbt]
    \centering
    \caption{Covering System for Sierp\'inski}
    \label{Scov}
    \begin{tabular}{|c|c|c|}
    \hline Covering & $p$ & $B\equiv \beta  \pmod p$\\
    \hline \hline
    - & $2^*$ & $1\pmod {2}$\\
     \hline
     $k\equiv 0 \pmod {2}$ & $3^*$ & $2\pmod {3}$\\
     \hline
     $k\equiv 3 \pmod {4}$ & $5^*$ & $3\pmod {5}$\\
     \hline
     $k\equiv 5 \pmod {8}$ & 17 & $9\pmod {17}$\\
     \hline
     $k\equiv 9 \pmod {16}$ & 257 & $129\pmod {257}$\\
     \hline
     $k\equiv 17 \pmod {32}$ & 65537 & $32769\pmod {65537}$\\
     \hline
     $k\equiv 1 \pmod {64}$ & 641 & $320\pmod {641}$\\
     \hline
     $k\equiv 33 \pmod {64}$ & 6700417 & $3350209\pmod {6700417}$\\
     \hline
    \end{tabular}
\end{table}
\end{center}

\begin{center}
\begin{table}[!hbt]
    \centering
    \caption{Covering System for Riesel}
    \label{Rcov}
    \begin{tabular}{|c|c|c|}
    \hline Covering & $p$ & $B \equiv \beta \pmod p$\\
    \hline \hline
    - & $2^*$ & $1\pmod {2}$\\
     \hline
     $k\equiv 1 \pmod {2}$ & $3^*$ & $2\pmod {3}$\\
     \hline
     $k\equiv 0 \pmod {3}$ & $7^*$ & $1\pmod {7}$\\
     \hline
     $k\equiv 0 \pmod {5}$ & 31 & $1\pmod {31}$\\
     \hline
     $k\equiv 2 \pmod {10}$ & $11^*$ & $3\pmod {11}$\\
     \hline
     $k\equiv 1 \pmod {15}$ & 151 & $76\pmod {151}$\\
     \hline
     $k\equiv 11 \pmod {45}$ & 631 & $57\pmod {631}$\\
     \hline
     $k\equiv 26 \pmod {45}$ & 23311 & $11446\pmod {23311}$\\
     \hline
     $k\equiv 86 \pmod {90}$ & 18837001 & $16\pmod {18837001}$\\
     \hline
     $k\equiv 8 \pmod {20}$ & 41 & $37\pmod {41}$\\
     \hline
     $k\equiv 18 \pmod {40}$ & 61681 & $61677\pmod {61681}$\\
     \hline
     $k\equiv 38 \pmod {80}$ & 4278255361 & $4278255357\pmod {4278255361}$\\
     \hline
     $k\equiv 78 \pmod {160}$ & 414721 & $414717\pmod {414721}$\\
     \hline
     $k\equiv 158 \pmod {160}$ & 44479210368001 & $4\pmod {44479210368001}$\\
     \hline
     $k\equiv 4 \pmod {30}$ & 331 & $269\pmod {331}$\\
     \hline
     $k\equiv 14 \pmod {60}$ & 61 & $39\pmod {61}$\\
     \hline
     $k\equiv 44 \pmod {60}$ & 1321 & $807\pmod {1321}$\\
     \hline
    \end{tabular}
\end{table}
\end{center}

\begin{thebibliography}{99} 

\bibitem{minmod} P.~Balister, B.~Bollob\'{a}s, R.~Morris, J.~Sahasrabudhe, and M.~Tiba, \textit{Erd\H{o}s covering systems}, Acta Mathematica Hungarica 161(1) (2020), 540--549.

\bibitem{barn}
G.~Barnes, author, ``Riesel conjectures and proofs", On page of No Prime Left Behind Project at http://www.noprimeleftbehind.net/crus/Riesel-conjectures.htm, Updated Sept. ~6, 2021.

\bibitem{harr1} C.~Bispels, M.~Cohen, J.~Harrington, J.~Lowrance, K.~Pontes, L.~Schaumann, and W.~H.~Wong, \textit{On Sierpiński and Riesel Repdigits and Repintegers}, Integers 26 (2026).


\bibitem{brucald} A.~Brunner and C.~Caldwell, D.~Krywaruczenko and C. ~Lownsdale, \textit{Generalized Sierpi\'{n}ski numbers base $b$}, preprint on University of Tennessee at Martin page (2008), https://www.utm.edu/staff/caldwell/preprints/2to100.pdf.

\bibitem{maria} M.~Cummings, M.~Filaseta, and O.~Trifonov, \textit{An upper bound for the minimum modulus in a covering system with squarefree moduli}, Acta Math. Hungarica 175 (2025), 1--25.

\bibitem{erdos} P.~Erd\H{o}s, \textit{On integers of the form $2^k+p$ and some related problems}, Summa Brasiliensis Mathematicae 2 (1950), 113-123.

\bibitem{fil}
M.~Filaseta, C.~Finch, and M.~Kozek,
\textit{On powers associated with Sierpi\'nski numbers, Riesel numbers and Polignac's conjecture}, 
J.~Number Theory 128 (2008), 1916--1940.

\bibitem{jacob}
M.~Filaseta and J.~Juillerat,
\textit{Consecutive primes which are widely digitally delicate}, INTEGERS: Ron Graham Memorial Volume, Vol.~21A, 2021, Paper No. A12, 37 pp.;
also see, Number Theory and Combinatorics: A Collection in Honor of the Mathematics of Ronald Graham, edited by Bruce M. Landman, Florian Luca, Melvyn B. Nathanson, Jaroslav Ne\v{s}et\v{r}il  and Aaron Robertson, Berlin, Boston: De Gruyter, 2022, 209--248.

\bibitem{Fildata} M.~Filaseta and J.~Juillerat, 
Data for ``Consecutive primes which are widely digitally delicate,"
{https://people.math.sc.edu/filaseta/Consecu-tiveWDDPrimes.html}. 

\bibitem{Fildatatwo} M.~Filaseta, J.~Juillerat, and T.~Luckner, 
Data for ``Consecutive primes which are widely digitally delicate and Brier numbers,"
{https://people.math.sc.edu/filaseta/ConsecutiveWDDBrierNumbers.html}. 

\bibitem{ftj} M.~Filaseta, J.~Juillerat, and T.~Luckner, \textit{Consecutive primes which are widely digitally delicate and Brier numbers}, Integers, 23 (2023).

\bibitem{jj}
M.~Filaseta, J.~Juillerat, and J.~Southwick, 
Widely Digitally Stable Numbers, in \textit{Combinatorial and Additive Number Theory IV} (ed.~M.~Nathanson), 
Springer Proc.~Math.~Stat.~347,
Springer, Cham, 2021, 161--193.

\bibitem{jerm} 
M.~Filaseta and J.~Southwick,
\textit{Primes that become composite after changing an arbitrary digit},
Math. Comp.~90 (2021), 979--993.

\bibitem{harr2} 
K.~Franzone and J.~Harrington,
\textit{Sierp\'inski and Riesel numbers in Narayana’s cow sequence},
Integers 24 (2024).

\bibitem{slobrier} 
O.~Gerard, author, The On-Line Encyclopedia of Integer Sequences, published electronically at https://oeis.org/A076335, Nov.~7, 2002.

\bibitem{hough}
B.~Hough,
\textit{Solution of the minimum modulus problem for covering systems},
Annals of Math. iss. 1, 181 (2015), 361--382.

\bibitem{data} T.~Luckner and R.~Philpott, 
Data for ``d-Translated Unit Sensitive Primes,"
{https://tluckner.github.io/DTUSBrierNumbers/}. 

\bibitem{maynard}
J.~Maynard,
\textit{Dense clusters of primes in subsets},
Compositio Math. 152 (2016), 1517--1554.

\bibitem{poli}
A.~de~Polignac, 
\textit{Six propositions arithmologiques déduites du crible d'Ératosthène}, Nouvelles annales de math., 1st ser., 8 (1849), 423--429.
 
\bibitem{ries}
H.~Riesel, 
\textit{N\r{a}gra stora primtal}, Elementa 39 (1956), 258--260.

\bibitem{rom}
N.~P.~Romanoff, \textit{\"Uber einige Sätze der additiven Zahlentheorie}, Math. Ann. 109 (1934), 668-–678.

\bibitem{shiu}
D.~K.~L.~Shiu, \textit{Strings of congruent primes}, J.~Lond.~Math.~Soc.~61 (2000), 359--373.

\bibitem{sier}
W.~Sierpi\'nski,
\textit{Sur un probl\`eme concernant les nombres $k\cdot 2^{n}+1$}, 
Elem.~Math.~15 (1960), 73--74.

\bibitem{slosier} 
N.~J.~A.~Sloane, author, The On-Line Encyclopedia of Integer Sequences, published electronically at https://oeis.org/A076336, Nov.~7, 2002.

\bibitem{slories} 
D.~W.~Wilson, editor, The On-Line Encyclopedia of Integer Sequences, published electronically at https://oeis.org/A101036, Jan.~17, 2005.

\end{thebibliography}
\end{document}